\documentclass[11pt,fleqn]{article}
\title{The Hilbert's class field  and the  $p$-class group of the cyclotomic fields}
\author{Roland Qu\^eme}
\usepackage{amsmath,amsthm,amstext,amsbsy,eufrak}
\newtheorem{thm}{Theorem}[section]

\newtheorem{lem}[thm]{Lemma}
\font\mathbb=msbm10
\newcommand{\N}{\mbox{\mathbb N}}

\newcommand{\Q}{\mbox{\mathbb Q}}
\newcommand{\Z}{\mbox{\mathbb Z}}

\newcommand{\modu}{\ \mbox{mod}\ }
\newcommand{\be}{\begin{equation}}
\newcommand{\ee}{\end{equation}}
\newcommand{\bd}{\begin{displaymath}}
\newcommand{\ed}{\end{displaymath}}
\newcommand{\bn}{\begin{enumerate}}
\newcommand{\en}{\end{enumerate}}

\newcommand{\mk}{\mathfrak}

\date{2011 jan 26}
\begin{document}
\maketitle
\tableofcontents
\abstract
Let $p$ be an irregular prime and  $K=\Q(\zeta)$ the $p$-cyclotomic field.
Let $\sigma$ be a $\Q$-isomorphism of $K$ generating $Gal(K/\Q)$.
Let $S/K$ be a cyclic unramified  extension of degree $p$, defined by  $S= K(A^{1/p})$ where $A\in K\backslash K^p$,  $A\Z_K=\mk a^p$ with  $\mk a$   non-principal  ideal of $\Z_K$, 
$A^{\sigma-\mu}\in K^p$
and $\mu\in{\bf F}_p$. We compute  explicitly the decomposition of the prime $p$ in the subfields $M$ of $S$ of degree $[M:\Q]=p$.

\section{Introduction}
Let $p$ be an irregular prime and  $K=\Q(\zeta)$ the $p$-cyclotomic field.
Let $\sigma$ be a $\Q$-isomorphism of $K$ generating $Gal(K/\Q)$.
Let $S/K$ be a cyclic unramified  extension of degree $p$, defined by  $S= K(\omega)$ where $\omega=A^{1/p}$,  $A$ is a primitive element of $K$,  $A\in K\backslash K^p$,  $A\Z_K=\mk a^p$ with  $\mk a$   non-principal  ideal of $\Z_K$, 
$A^{\sigma-\mu}\in K^p$ and $\mu\in{\bf F}_p$. 
The extensions $S/K$  and $K/\Q$ are  respectively  cyclic of degree $p$ and $p-1$. At each $p$th root $\omega$ of $A$ correspond a subfield $M$ of $S$ with  degree $[M:\Q]=p$, $S/M$ cyclic  and $M/\Q$ non galois.  In this article,  we  describe completely the decomposition of the prime $p$ in the  fields  extensions $(S/K, K/\Q)$  and  $(S/M, M/\Q)$.

%%% ====================================================================
%
%RRRRRRRR 20
%%% ====================================================================
\section{On Kummer theory}
In  this section we recall some definitions and classical properties of the $p$-cyclotomic field.
\subsection{Some definitions and notations}\label{s601191}
In this subsection   we introduce  some  definitions and  notations on
cyclotomic fields,  $p$-class group, singular  numbers,  primary and
non-primary, used in this note.
\begin{enumerate}
\item
Let $p$ be an odd irregular  prime. Let $\zeta$ be a  $p$th primitive root of
the unity. Let $K$ be the $p$-cyclotomic field $K=_{def}\Q(\zeta)$
and $\Z_K$ its ring of integers. Let $K^+$ be the maximal totally
real subfield of $K$, $\Z_{K^+}$ its ring of integers and
$\Z_{K^+}^*$ the group of unit of $\Z_{K^+}$. Let $v$ be a primitive
root $\modu p$ and $\sigma:_{def} \zeta\rightarrow \zeta^v$ be a
$\Q$-automorphism of $K$. Let $G$ be the Galois group of the
extension $K/\Q$. Let ${\bf F}_p$ be the finite field of cardinal
$p$ and   ${\bf F}_p^*$ its multiplicative group. Let
$\lambda=\zeta-1$ (By opposite in this note to the usual
notation $\lambda=\zeta-1$). The prime ideal of $K$ lying over $p$
is $\pi=\lambda O_K$. For $a\in K$ we adopt the notation  $\overline{a}$ for the 
complex conjugate of $a$.
\item
Let $C_p$ be the subgroup of exponent $p$ of the of the class group of $K$ (called 
$p$-class group in this article).

$C_p$ is the direct sum of $r$ subgroups  $\Gamma_i$  of order $p$,
each $\Gamma_i$ annihilated by a polynomial $\sigma-\mu_i\in {\bf
F}_p[G]$ with $\mu_i\in{\bf F}_p^*$,
\be\label{e611011}
C_p=\oplus_{i=1}^r \Gamma_i,
\ee
where $\mu\equiv v^{n}\modu p$ with
a natural integer $n,\quad 2\leq n\leq p-2$.
\item
$r_p$ is the rank of $C_p$ seen as a ${\bf F}_p[G]$-module. If $r_p>0$ then $p$ is irregular.
Let $h^+$ be the class number of $K^+$.
Let
$C_p^+$ be the $p$-class group of $K_p^+$. Let $r_p^+$ be the rank of
$C_p^+$. Then $C_p=C_p^+\oplus C_p^-$ where $C_p^-$ is the relative
$p$-class group. $r_p^+\geq 1$ if and only if $h^+\equiv 0\mod p$.
\item
A number  $A\in \Z_K$ is said singular if $v_\pi(A)=0$, $A^{1/p}\not\in K$
and if  there exists an   ideal $\mk a$  of $K$ such that
$A \Z_K =\mk a^p$. Observe that, with this definition, a unit $\eta\in \Z_{K^+}^*$
with $\eta^{1/p}\not \in \Z_{K^+}^*$ is singular.
\item
A number $A\in K$ is said {\it semi-primary} if $v_\pi(A)=0$ and if there exists a natural integer $a$ such that
$A\equiv a\modu\pi^2$.
A number $A\in K$ is said {\it primary} if $v_\pi(A)=0$ and if there exists a natural integer $a$ such that
$A\equiv a^p\modu\pi^p$. Clearly a primary number is semi-primary.
A number $A\in K$ is said {\it hyper-primary} if $v_\pi(A)=0$ and if there exists a natural integer $a$ such that
$A\equiv a^p\modu\pi^{p+1}$. \en
%%% ====================================================================
%
%RRRRRRRR 20
%%% ====================================================================
\subsection{Some preliminary results}
In this section we recall some classical properties of singular numbers (given for instance in Qu\^eme \cite{que}
in theorems  2.4 p. 4, 2.7 p. 7 and 3.1 p. 9). Let $\Gamma$ be one of the $r$ subgroups $\Gamma_i$ of order $p$  of $C_p$ defined in relation
(\ref{e611011}).
\begin{enumerate}
\item\label{i609092}
\underline {If  $r_p^->0$ and  $\Gamma\subset  C_p^-$: }
then there exist  singular semi-primary integers $A$  with $A  \Z_K =\mk  a^p$
where $\mk a$ is a {\bf non}-principal  ideal of $\Z_K$  verifying   simultaneously
\begin{equation}\label{e512101}
\begin{split}
& Cl(\mk  a)\in \Gamma,\ Cl(\mk  a^{\sigma-\mu})=1,\\
& \sigma(A)=A^\mu\times\alpha^p,\quad \mu\in {\bf F}_p^*,\quad \alpha\in K,\\
&\mu\equiv v^{2m+1}\modu p, \quad m\in\N, \quad 1\leq m\leq \frac{p-3}{2},\\
&\pi^{2m+1} \ |\ A-a^p,\quad a\in\N,\quad 1\leq a\leq p-1.\\
\end{split}
\end{equation}
In that case we say that $A$ is a {\it negative} singular integer to point out the fact that $Cl(\mk a)\in C_p^-$.
Moreover, this number $A$ verifies
\begin{equation}\label{e512103}
A\times\overline{A}=D^p,
\end{equation}
for some integer $D\in \Z_{K^+}^*$.
\begin{enumerate}
\item
Either  $A$ is singular non-primary  with  $\pi^{2m+1} \ \|\ A-a^p$.
\item
Or  $A$ is singular  primary  with   $\pi^{p}\   |\ A-a^p$. In that case we know by reflection from class field theory  that $r_p^+>0$.
\end{enumerate}
(see  for instance Qu\^eme \cite{que} theorem 2.4 p. 4).
\item\label{i609091}
\underline {If $r_p^+>0$ and    $\Gamma\subset C_p^+$: }
then there exist  singular semi-primary  integers $A$  with $A \Z_K =\mk a^p$
where $\mk a$ is a {\bf non}-principal  ideal of $\Z_K$ and verifying   simultaneously
\begin{equation}\label{e6012210}
\begin{split}
& Cl(\mk  a)\in \Gamma,\ Cl(\mk  a^{\sigma-\mu})=1,\\
& \sigma(A)=A^\mu\times\alpha^p,\quad \mu\in {\bf F}_p^*,\quad \alpha\in K,\\
&\mu\equiv v^{2m}\modu p, \quad m\in\Z, \quad 1\leq m\leq \frac{p-3}{2},\\
&\pi^{2m} \ |\ A-a^p,\quad a\in\Z,\quad 1\leq a\leq p-1,\\
\end{split}
\end{equation}
In that case we say that $A$ is a {\it positive } singular integer to point out the fact that $Cl(\mk  a)\in C_p^+$.
Moreover, this integer $A$ verifies
\begin{equation}\label{e512103}
\frac{A}{\overline{A}}=D^p,
\end{equation}
for some number  $D\in K$.  If $h^+\equiv 0\modu p$, it is important  to note   that $D\not=1,\  D\in K\backslash K^+$ is possible, for instance when
$\mk  a=\mk  q$, with  $\mk  q$  prime ideal of $\Z_K$,  $Cl(\mk  q)\in C_p^+$ and $q\equiv 1\modu p$ (from Kummer, the prime $q\equiv 1 \mod p$ generates the class group of $K$).
\begin{enumerate}
\item
Either  $A$ is singular non-primary with  $\pi^{2m} \ \|\ A-a^p$ .
\item
Or  $A$ is singular  primary with  $\pi^{p}\  |\ A-a^p$.
\end{enumerate}
(see  for instance  \cite{que} theorem 2.7 p. 7).
\en

%%% ====================================================================
%
%RRRRRRRR 20
%%% ====================================================================
\section{Singular extensions}
%\subsection {Some Definitions}
We introduce a terminology  used in this note:
\bn
\item
In this section, let us fix   $\Gamma$  one of the $r$ subgroups of order $p$ of $C_p$ defined by relation (\ref{e611011}).
Let $A$ be a singular semi-primary  integer,    negative or  positive,  verifying   respectively  the relations
(\ref{e512101}) or (\ref{e6012210}). Without loss of generality we assume in the sequel that $A\equiv 1\mod\pi^2$ because  if $A\not\equiv 1\mod p^2$ we can replace 
$A$ by $A^\prime= A^{p-1}$.
Let us define  $$S=_{def} \Q(A^{1/p}).$$
The singular integer  $A\in K$, so the absolute degree $[S:\Q]$ is $p\times deg(A)$ dividing $p(p-1)$ where $deg(A)$ is the degree of the minimal polynomial of the algebraic integer $A$.
{\it In the sequel, we consider  the only cases where $A$ is of degree $p-1$ (or a primitive element of $K$).}

We call $S/K$ a  singular   extension. The singular extension $S/K$ is negative or positive depending on whether $A$ is a negative or positive singular number.
From Kummer theory $S/K$ is a cyclic extension of degree $p$ and $S/\Q$ is a Galois non-abelian extension of degree $(p-1)p)$.
\item
A singular extension $S=K(A^{1/p})$ is said primary or  non-primary
depending on whether the singular number $A$  is primary or non-primary.
From Hilbert class field theory, if $S$ is  primary,  the extension  $S/K$ is
the cyclic unramified extension of degree $p$ corresponding to the cyclic group $\Gamma$ defined in relation (\ref{e611011})
 p. \pageref{e611011}. 

Observe that if $A$ is singular primary then $S/K$ is unramified and the principal ideal $\pi$ of $K$ splits totally in $S/K$ from
principal ideal theorem which implies that  $A$ is singular hyperprimary.
\en

%%% ====================================================================
%
%RRRRRRRR 20
%%% ====================================================================

\section{Singular $\Q$-fields }
Let $A$ be a semi-primary  singular  integer, negative (see definition (\ref{e512101})), positive (see definition (\ref{e6012210})).
Let $\omega$ be a $p$th root of $A$
\be\label{e610096}
\omega:= A^{1/p}.
\ee
Then $S/K=K(\omega)/K$ is  the corresponding singular $K$-extension.
Observe that this definition implies that
$\omega\in \Z_{S}$ ring of integers of $S$.
%%% ====================================================================
%
%RRRRRRRR 20
%%% ====================================================================
%%% ====================================================================
%
%RRRRRRRR 20
%%% ====================================================================
\begin{lem}\label{l609101}
Suppose that $S/K$ is a singular primary extension.
Let $\theta :\omega\rightarrow\omega\zeta$ be a $K$-isomorphism of the field $S$.
There are $p$ prime ideals  of $\Z_{S}$ lying over $\pi$.
There  exists a prime ideal $\pi_0$ of $\Z_{S}$ lying over $\pi$ such that
the $p$ prime ideals $\pi_n=\theta^n(\pi_0),\ n=0,\dots,p-1$ of $\Z_{S}$ lying over $\pi$ verify the congruences
\be\label{e609043}
 \begin{split}
 &\pi_0^ 2\ |\ \omega-1,\\
 &\pi_n\ \|\ \omega-1, \ \dots, n=1,\dots,p-1.\\
\end{split}
 \ee
  \begin{proof}$ $
\bn
\item
From Hilbert class field theory and Principal Ideal Theorem, the prime principal ideal $\pi$ of $K$ splits totally in the  unramified cyclic extension $S/K$.
It follows that $\pi^{p+1}\|A-1$ (see for instance Ribenboim \cite{rib} case III p.168).
\item
Then $\omega^p-1\equiv 0\modu\pi^{p+1}$.
Let $\theta :\omega\rightarrow \omega\zeta$ be a $K$-automorphism of the field  $S$.
Let $\Pi$ be any  of the $p$ prime ideals of $\Z_{S}$ lying over $\pi$.
Then  $\pi \Z_{S}=\prod_{l=0}^{p-1}\theta^l(\Pi)$ where
$\theta^l(\Pi),\ l=0,\dots,p-1$ are the $p$ prime ideals of $\Z_{S}$ lying over $\pi$.
\item
From $A\equiv 1\modu \pi^{p+1}$ we see that
\bd
\omega^p-1=\prod_{i=0}^{p-1}(\omega\zeta^{i}-1)\equiv 0\modu\prod_{l=0}^{p-1}\theta^l(\Pi)^{p+1}.
\ed
It follows that there exists at least one $i$ such that $\omega\zeta^{i}-1\equiv 0\mod \Pi^2$.
There exists only one such $i$  if not, from $\omega\zeta^{i^\prime}-1\equiv 0\mod \Pi^2$, we should have
$\zeta^i-\zeta^{i^\prime}\equiv 0\mod\Pi^2$.
Therefore $\omega-1\equiv 0\mod \theta^{-i}(\Pi)$. The theorem follows with $\pi_0=\theta^{-i}(\Pi)$.
 \en
 \end{proof}
 \end{lem}
%% ====================================================================
%
%RRRRRRRR 20
%%% ====================================================================
%%% ====================================================================
%
%RRRRRRRR 20
%%% ====================================================================
\begin{lem}\label{l609121}
Suppose that $S/K$ is a singular non primary  extension. Then $\pi$ ramifies in $S/K$ and  $\Pi\ |\ \omega-1$ where $\Pi$ is  the prime of $S$ lying over $\pi$. 
\begin{proof}$ $
\bn
\item
Show that the extension $S/K$ is ramified at $\pi$: $S/K$ is unramified at all the primes except $p$ because $A$ is singular. $\pi$ does not splits 
in $S/K$, if not $S/K$ should be a singular primary extension. $\pi$ is not inert in $S/K$,  if not we should have $\omega\equiv 1\mod \pi$
so $A\equiv 1\mod \pi^p$ and $S/K$ should be primary, therefore $\pi$ ramifies in $S/K$ and  $\pi \Z_{S}=\Pi^p$.
\item
$A\equiv 1\modu\pi$,  thus $\omega\equiv 1\modu\Pi$.
\en
\end{proof}
\end{lem}
%% ====================================================================
%
%RRRRRRRR 20
%%% ====================================================================
There are $p$ different automorphisms   of the field $S$ extending the $\Q$-automorphism $\sigma$ of the field $K$, we will fix one of them.
\begin{lem}\label{l609095}
There exists  one and only one automorphism  $\sigma_\mu$ of  $S/\Q$ extending $\sigma$  such that
\be\label{e609095}
\omega^{\sigma_\mu-\mu}\equiv 1\modu \pi^2.
\ee
\begin{proof}
\item
From $\sigma(A)=A^\mu\alpha^p$ there exist $p$ different automorphisms $\sigma_{(w)},\ w=0,\dots,p-1,$ of the field  $S$
extending the $\Q$-isomorphism
$\sigma$ of the field $K$,  defined by
\be\label{e610137}
\sigma_{ (w)}(\omega)=\omega^\mu\alpha\zeta^w,
\ee
for  natural numbers $w=0,1,\dots,p-1$.
There exists one and only one  $w$ such that $\alpha\times \zeta^w$ is a semi-primary number (or
$\alpha\times\zeta^w\equiv 1\modu\pi^2$). Let us set $\sigma_\mu=\sigma_{(w)}$ to emphasize the role of $\mu$.
Therefore we get
\be\label{e609041}
\sigma_\mu(\omega)\equiv \omega^\mu\modu\pi^2.
\ee
\end{proof}
\end{lem}
In the sequel, without loss of generality, we fix the semi-primary number $\alpha\in K$ such that
\be\label{e809041}
\sigma_\mu(\omega)=\omega^\mu\alpha.
\ee
%%% ====================================================================
%
%RRRRRRRR 20
%%% ====================================================================
\begin{lem}\label{l609096}
$\sigma_\mu^{p-1}(\omega)=\omega$.
\begin{proof}
We have $\sigma_\mu^{p-1}(A)= A$, therefore there exists a natural integer $w_1$ such that
$\sigma_\mu^{p-1}(\omega)=\omega\times \zeta^{w_1}$.
We have proved in relation (\ref{e609041}) that
\be\label{e609103}
\sigma_\mu(\omega)\equiv \omega^\mu \modu \pi^2,
\ee
thus $\sigma_\mu^{p-1}(\omega)\equiv\omega^{\mu^{p-1}}\equiv\omega\times  A^{(\mu^{p-1}-1)/p}\equiv \omega\modu \pi^2$
which implies that $w_1=0$ because $A$ is semi-primary. Therefore  $\sigma_\mu^{p-1}(\omega)=\omega$.
\end{proof}
\end{lem}
%%% ====================================================================
%
%RRRRRRRR 20
%%% ====================================================================
Let us define $\Omega\in \Z_{S}$ ring of integers of $S$ by the relation
\be\label{e609061}
\Omega=\sum_{i=0}^{p-2}\sigma_\mu^i(\omega).
\ee
%%% ====================================================================
%
%RRRRRRRR 20
%%% ====================================================================
\begin{thm}\label{l609097}
$M=\Q(\Omega)$ is a field with $[M:\Q]= p$, $[S: M]=p-1$  and $\sigma_\mu(\Omega)=\Omega$.
\begin{proof}$ $
\bn
\item
Show that $\Omega\not=0$: If $S/K$ is unramified, then $\omega\equiv 1\modu\pi$ implies with definition of $\Omega$ that
$\Omega\equiv p-1\modu \pi$ and so $\Omega\not=0$. If $S/K$ is ramified, then $\omega\equiv 1\modu\Pi$
implies with definition of $\Omega$ that
$\Omega\equiv p-1\modu \Pi$ because $\sigma_\mu(\Pi)=\Pi$ and so $\Omega\not=0$.
\item
Show that $\Omega\not\in K$: from $\sigma_\mu(\omega)=\omega^\mu\alpha$ we get
\bd
\Omega=\sum_{i=0}^{p-2}\omega^{\mu^i\mod p}\times \beta_i,
\ed
with $\beta_i\in K$.
Putting together terms of same degree   we get  $\Omega=\sum_{j=1}^{p-1}\gamma_j\omega^j$ where $\gamma_j\in K$ are not all null
because $\Omega\not=0$.
$\Omega\in K$
should imply the polynomial equation $\sum_{j=1}^{p-1}\omega^{j}\times\gamma_j-\gamma=0$ with $\gamma\in K$,  not possible
because the minimal polynomial equation of $\omega$ with coefficients in $K$ is $\omega^p-A^{p-1}=0$.
\item
Show that $M=\Q(\Omega)$ verifies $M\subset S$ with $[M:\Q]=p$ and  $[S:M]=p-1$:
$S/\Q$ is a Galois extension with $[S:\Q]=(p-1)p$.
Let $G_S$ be the Galois group of $S/\Q$. Let $<\sigma_\mu>$ be the subgroup of $G_S$
generated by the automorphism $\sigma_\mu\in G_S$. We have seen in lemma \ref{l609096} that $\sigma_\mu^{p-1}(\omega)=\omega$.
In the other hand $\sigma_\mu^{p-1}(\zeta)=\zeta$ and $\sigma_\mu^n(\zeta)\not=\zeta$
for $n<p-1$ and so $<\sigma_\mu>$ is of order $p-1$.
\item
From fundamental theorem of Galois theory,  there is a fixed field $M=S^{<\sigma_\mu>}$  with
$[M:\Q]= [G_S:<\sigma_\mu>]=p$.
From $\sigma_\mu(\Omega)=\Omega$ seen and from definition relation (\ref{e609061}) it follows that $\Omega\in M$ and from $\Omega \not\in K$
it follows that $M=\Q(\Omega)$.
Thus  $S=M(\zeta)$
and    $\omega\in S$ can be written
\be\label{e608093}
\omega=\sum_{i=0}^{p-2} \omega_i\lambda^i,\ \omega_i\in M.
\ee
with $\lambda=\zeta-1$ and with  $\sigma_\mu(\omega_i)=\omega_i$  because $\sigma_\mu(\Omega)=\Omega$.
\en
\end{proof}
\end{thm}
%%% ====================================================================
%
%RRRRRRRR 20
%%% ====================================================================
\paragraph{Some  definitions:}
The field  $M\subset S$ is called a singular $\Q$-field.  A singular $\Q$-field $M$ is said primary (respectively  non-primary)
if $S$ is a singular primary (respectively  non-primary) extension.

%%% ====================================================================
%
%RRRRRRRR 20
%%% ====================================================================
\section{The decomposition of the prime $p$  in the singular primary $\Q$-fields}

This  section deals with the decomposition of the prime $p$ in the  singular primary  $\Q$-fields  $M$.
In that case  recall that $S/K$ is a cyclic unramified extension of degree $p$ and there
are $p$ prime ideals in $S/K$ over $\pi$ 

Observe  that the case of singular non-primary $\Q$-fields can easily be described. The extension $S/K$ is fully ramified at $\pi$,
so $p \Z_{S}=\pi_{S}^{p(p-1)}$. Therefore there is only one prime ideal $\mk p$ of $M$ ramified
with $p \Z_{M}=\mk p^p$.

%%% ====================================================================
%
%RRRRRRRR 20
%%% ====================================================================
Recall that $\theta$ is the $K$-isomorphism $\theta: \omega\rightarrow \omega\zeta$ of $S$ and  that $\pi$ splits totally in $S/K$ with
$\pi_i=\theta^i(\pi_0)$ for $i=0,\dots,p-1$.
\begin{lem}\label{l609102}
$\sigma_\mu(\pi_0)=\pi_0$
\begin{proof}
From relation (\ref{e609103}) $\sigma_\mu(\omega)\equiv\omega^\mu\modu\pi^2$. From lemma \ref{l609101}, $\omega\equiv 1\modu\pi_0^2$
and so $\sigma_\mu(\omega)\equiv\omega^\mu\equiv 1\modu\pi_0^2$.  Then $\omega\equiv 1\modu \sigma_\mu^{-1}(\pi_0)^2$.
If $\sigma_\mu^{-1}(\pi_0)\not=\pi_0$ it follows that $\omega\equiv 1\modu\pi_0^2\times\sigma^{-1}(\pi_0^2)$, which contradicts lemma
\ref{l609101}.
\end{proof}
\end{lem}
%%% ====================================================================
%
%RRRRRRRR 260
%%% ====================================================================
\begin{lem}\label{l301281}
Let $\pi_k=\theta^k(\pi_0)$ for any $k\in\N,\quad 1\leq k\leq p-1$.
Then $\sigma_\mu(\pi_k)=\pi_{n_k}$ with $n_k\in\N, \quad n_k\equiv k\times v \mu^{-1}\modu p$.
\begin{proof}$ $
\bn
\item
From $\pi_0^{ 2}\  |\ \omega-1$, it follows that $$\theta^k(\pi_0^{ 2})=\pi_k^{ 2}\ |\ \omega\zeta^k-1.$$
Then
\bd
\sigma_\mu(\pi_k)^{ 2}\ |\ \sigma_\mu(\omega)\times\zeta^{v k}-1.
\ed
\item
We have $\sigma_\mu(\pi_k )=\pi_{k+l_k}$ for some $l_k\in\N$ depending on $k$.
From relation (\ref{e609103}) we know that $\sigma_\mu(\omega)\equiv \omega^\mu \modu\pi^2$.
Therefore
\begin{displaymath}
\pi_{k+l_k}^{ 2}\ |\ \omega^\mu\times\zeta^{v k}-1.
\end{displaymath}
\item
In an other part by the $K$-automorphism $\theta^{k+l_k}$ of $S$ we have
\bd
\pi_{k+l_k}^{ 2}\ |\ \omega\times\zeta^{k+l_k}-1,
\ed
so
\begin{displaymath}
\pi_{k+l_k}^{ 2}\ |\ (\omega^\mu\times\zeta^{\mu(k+l_k)}-1).
\end{displaymath}
\item
Therefore
$\pi_{k+l_k}^{ 2}\ |\  \omega^\mu(\zeta^{v k}-\zeta^{\mu(k+l_k)})$,
and so
\begin{displaymath}
\pi_{k+l_k}^{ 2}\ |\ \zeta^{v k}-\zeta^{\mu(k+l_k)}.
\end{displaymath}
\item
This implies that  $\mu(k+l_k)-v k\equiv 0\modu p$, so
$\mu l_k+k(\mu-v)\equiv 0\modu p$ and
finally that
\begin{displaymath}
l_k\equiv k\times\frac{v-\mu}{\mu}\mod p,
\end{displaymath}
where we know that $v-\mu\not\equiv 0\modu p$ from Stickelberger relation. Then
$n_k\equiv k+k\times\frac{v-\mu}{\mu}=k\times\frac{v}{\mu}\modu p$, which achieves the proof.
\en
\end{proof}
\end{lem}
%%% ====================================================================
%
%RRRRRRRR 260
%%% ====================================================================
\begin{lem}\label{l609271}$ $
\bn
\item
If $S/K$ is a  singular primary negative extension then $\sigma_\mu^{(p-1)/2}(\pi_k)=\pi_k$.
\item
If $S/K$ is a singular primary positive   extension then $\sigma_\mu^{(p-1)/2}(\pi_k)=\pi_{n-k}$.
\en
\begin{proof}
From lemma \ref{l301281} we have $\sigma_\mu^{(p-1)/2}(\pi_k)= \pi_{k^\prime}$ with
$k^\prime\equiv k v^{(p-1)/2}\mu^{-(p-1)/2}\mod p$.
 If $S/K$ is negative then  $v^{(p-1)/2}\mu^{-(p-1)/2}\equiv 1\modu p$ and if
$S/K$ is positive or unit then  $v^{(p-1)/2}\mu^{-(p-1)/2}\equiv -1\modu p$
and the result follows.
\end{proof}
\end{lem}
%%% ====================================================================
%
%RRRRRRRR 260
%%% ====================================================================
\begin{lem}\label{l609105}
The length of the orbit of the action of the group $<\sigma_\mu>$ on $\pi_0$ is $1$
and the length of the orbit of the action of the group $<\sigma_\mu>$ on
$\pi_i, \  i=1,\dots,p-1$ is $d$ where $d$ is the order of $v\mu^{-1}\modu p$.
\begin{proof}
For $\pi_0$ see lemma \ref{l609102}.
For $\pi_k$ see lemma \ref{l301281}: $\sigma_\mu(\pi_k)=\sigma_\mu(\pi_{n_k})$ with
$n_k\equiv v\mu^{-1}\modu p$, then $\sigma_\mu^2(\pi_k)=\sigma_\mu(\pi_{n_{k_2}})$ with
$n_{k_2}\equiv k v^2\mu^{-2}\modu p$ and finally $n_{k_d}\equiv k \modu p$.
\end{proof}
\end{lem}
%%% ====================================================================
%
%RRRRRRRR 260
%%% ====================================================================
The only prime ideals of $M/\Q$ ramified are lying over $p$.
The prime ideal of $K$ over p is $\pi$.
To avoid cumbersome notations,
the prime ideals of $S$ over $\pi$ are noted here
$\Pi$ or $\Pi_i=\theta^i(\Pi_0), \ i=1,\dots,p-1$,
and the prime ideals of $M$ over $p$ are noted $\mk P$ or  $\mk P_j,\ j=0,\dots,\nu$ where   $\nu+1$ is the number of such ideals.
%%% ====================================================================
%
%RRRRRRRR 260
%%% ====================================================================
\begin{thm}\label{t609101}$ $
Let  $d$ be  the order of $v\mu^{-1}\modu p$.
There are $\frac{p-1}{d}+1$ prime ideals in the singular primary $\Q$-field $M$ lying over $p$.
Their  prime decomposition and ramification is:
\bn
\item
$ e(\mk P_0/p\Z) =1.$
\item
$e(\mk P_j/p\Z)=d$ for all $j=1,\dots,\frac{p-1}{d}$  with $d>1$.
\en
\begin{proof}$ $
\bn
\item
\underline{preparation of the proof}
\bn
\item
The inertial degrees verifies
$f(\pi/p\Z)=1$ and $f(\Pi/\pi)=1$ and so $f(\Pi/p\Z)=1$.
Therefore, from multiplicativity of degrees in extensions,
it follows that
$f(\mk P/p\Z)=f(\Pi/\mk P)=1$
where $\Pi$ is lying over $\mk P$.
\item
$e(\pi/p\Z)=p-1$ and $e(\Pi/\pi)=1$ and so $e(\Pi/p\Z)=p-1$.
\item
Classically, we get
\be\label{e609274}
\sum_{j=0}^\nu e(\mk P_j/p\Z)=p,
\ee
where $\nu+1$ is the number of prime ideals of $M$ lying over $p$ and
where $e(\mk P_j/p\Z)$ are ramification indices dividing $p-1$ because, from multiplicativity of degrees in extensions,
$e(\Pi/\mk P_j)\times e(\mk P_j/p\Z) =p-1$.
\en
\item
\underline{Proof}
\bn
\item
The extension $S/M$ is Galois of degree $p-1$, therefore
the number $n_\mk P$ of prime ideals $\Pi$ of $\Z_S$ lying over one $\mk P$ of $\Z_M$ is $n_\mk P=\frac{p-1}{e(\Pi/\mk P)}=e(\mk P/p\Z)$.
\item
Let $c(\Pi)$ be the orbit of $\Pi$ under the action of the group $<\sigma_\mu>$ of cardinal $p-1$ seen in lemma \ref{l609105}.
If $\Pi=\pi_0$ then the orbit $C_\Pi$ is of length $1$.
If $\Pi\not=\pi_0$ then the orbit $C_\Pi$ is of length $d$.
If $C_\Pi$ has one ideal lying over $\mk P$ then it has all its $d$ ideals lying over $\mk P$ because $\sigma_\mu(\mk P)=\mk P$.
This can be extended to   all $\Pi^\prime$ lying over $\mk P$ with $C_{\Pi^\prime}\not= C_\Pi$  and it follows that
when $\Pi\not=\pi_0$ then 
\be\label{e1101231}
d\ |\  n_\mk P=e(\mk P/p\Z).
\ee
There is  one $\mk P$ with $n_\mk P=e(\mk P/p\Z)=1$ because $C_{\pi_0}$ is the only orbit with one element.
\item
The extension $S/M$ is cyclic of degree $p-1$. There exists one field $N$ with
$M\subset N\subset S$ with degree  $[N:M]=\frac{p-1}{d}$.
If there were   at least two different prime  ideals $\mk P_{1}^\prime $ and $\mk P_{2}^\prime$ of $N$ lying over $\mk P$
it should follow
that $\mk P_{2}^\prime=\sigma_\mu^{j d}(\mk P_{1}^\prime)$ for some $j,\ 1\leq j\leq \frac{p-1}{d}$,
because the Galois group of $S/M$ is $<\sigma_\mu>$ and the Galois group of $N/M$ is $<\sigma_\mu^{d}>$.
But, if a prime ideal $\pi_k$ of $S$ lies over $\mk P_{1}^\prime$ then $\sigma_\mu^{jd}(\pi_k)$ should lie over
$\mk  P_{2}^\prime$.
From lemma \ref{l609105},  $ \sigma_\mu^{d}(\pi_k)=\pi_k$ should  imply that $\mk P_{2}^\prime=\mk P_{1}^\prime$, contradiction.
Therefore the only possibility is that $\mk P$ is fully ramified in $N/M$ and
thus  $\mk P O_{N}=\mk P^{\prime (p-1)/d}$.
Therefore $e(\mk P^\prime/\mk P)=\frac{p-1}{d}$ and so
$e(\mk P^\prime/p\Z)=e(\mk P/p\Z)\times \frac{p-1}{d}\ |\ p-1$ and thus $$e(\mk P/p\Z)\ |\ d.$$
From previous result (\ref{e1101231})  it follows that $e(\mk P/p\Z)=d$. Then $d>1$ because $\mu-v\not\equiv 0\modu p$ from Stickelberger theorem.
There are $\frac{p-1}{d}+1$  prime ideals $\mk P_i$ because, from relation (\ref{e609274})
$p= 1+\sum_{i=1}^\nu e(\mk P_i/p\Z)=1+\nu \times d$.
\en
\en
\end{proof}
\end{thm}
%% ====================================================================
%
%RRRRRRRR 20
%%% ====================================================================
\paragraph{Example:}
let us consider the case of prime numbers $p$ with  $\frac{p-1}{2}$  prime.
\bn
\item\label{i1012151}
\underline{Singular primary negative  $\Q$-fields}

Here $\mu^{(p-1)/2}\equiv -1\modu p$ and
$d\in\{2,\frac{p-1}{2},p-1\}$.
Straightforwardly,  the case $d=2$ is not possible: $\mu^2\equiv v^2\modu p$, then $\mu+v\equiv 0\modu p$ because
$\mu\not\equiv v\modu p$, then $\mu^{(p-1)/2}+v^{(p-1)/2}\equiv 0\modu p$, contradiction
because $\mu^{(p-1)/2}=v^{(p-1)/2}=-1$. $d=p-1$ is not possible because  $\mu^{(p-1)/2}-v^{(p-1)/2}\equiv 0\modu p$.
Therefore $d= \frac{p-1}{2}$, so  the ramification of $p$ in the singular $\Q$-field $M$ is
$e(\mk P_0/p\Z)=1$ and $e(\mk P_1/p\Z)=e(\mk P_2/p\Z)=\frac{p-1}{2}$.
\item
\underline{Singular primary positive   $\Q$-fields}

Here $\mu^{(p-1)/2}\equiv 1\modu p$ and
$d\in\{2,\frac{p-1}{2},p-1\}$.
The case $d=2$  : $\mu^2-v^2\equiv 0\modu p$ then $\mu+v\equiv 0\modu p$ so $\mu\equiv v^{(p+1)/2}\modu p$,
so $B_{(p+1)/2}\equiv 0\modu p$
where $B_{(p+1)/2}$ is a Bernoulli Number.
The case $d=\frac{p-1}{2}$ is not possible because $\mu^{(p-1)/2}\equiv 1\modu p$ and $v^{(p-1)/2}\equiv -1\modu p$.
The case  $d=p-1$ is possible, so the ramification of $p$ in the singular $\Q$-field $M$ is
\bn
\item either $e(\mk P_0/p\Z)=1$
and $e(\mk P_i/p\Z)=2$ for $i=1,\dots,\frac{p-1}{2}$ and $B_{(p+1)/2}\equiv 0\mod p$.
\item or $ e(\mk P_0/p\Z)=1$ and $e(\mk P_1/p\Z)=p-1$.
\en
\en

%%% ====================================================================

%% ====================================================================
%
%RRRRRRRR 20
%%% ====================================================================
Roland Qu\^eme

13 avenue du ch\^ateau d'eau

31490 Brax

France

mailto: roland.queme@wanadoo.fr

************************************

MSC Classification : 11R18;  11R29; 11R32

************************************
%% ====================================================================
%
%RRRRRRRR 20
%%% ====================================================================
\end{document}